\documentclass[12pt]{amsart}
\usepackage[latin1]{inputenc}
\usepackage[english]{babel}
\usepackage{amsmath,amssymb,amsthm,amsfonts, color}
\usepackage[all]{xy}
\usepackage{tikz}
\usepackage{tikz-cd}

\usepackage{latexsym,dsfont}
\usepackage{enumerate}
\usepackage{graphicx}
\usepackage{latexsym}

\makeatletter
\@namedef{subjclassname@2010}{ \textup{2010} Mathematics Subject Classification}
\makeatother

\newtheorem{theorem}{Theorem}

\newtheorem{lemma}{Lemma}

\theoremstyle{definition}

\newtheorem{remark}[lemma]{Remark}

\numberwithin{equation}{section}

\numberwithin{equation}{section}

\DeclareMathOperator{\esssup}{ess~sup}

\begin{document}
\title{Abstract Ces\`aro spaces. II. Optimal range{\rm *}}\thanks{{\rm *}This publication has been written during scholarship 
period of the first author at the Lule{\aa} University of Technology, thanks to a Swedish Institute scholarschip 
(number 0095/2013).}
\begin{abstract}
Abstract Ces\`aro spaces are investigated from the optimal domain and optimal range point of view. There is 
a big difference between the cases on $[0, \infty)$ and on $[0, 1]$, as we can see in Theorem 1. Moreover, we present an improvement of Hardy inequality on  $[0, 1]$ which plays an important role in these considerations. 
\end{abstract}
\author[Le\'snik]{Karol Le\'snik}
\address[Karol Le{\'s}nik]{Institute of Mathematics\\
of Electric Faculty Pozna\'n University of Technology, ul. Piotrowo 3a, 60-965 Pozna{\'n}, Poland}
\email{\texttt{klesnik@vp.pl}}
\author[Maligranda]{Lech Maligranda}
\address[Lech Maligranda]{Department of Engineering Sciences and Mathematics\\
Lule{\aa} University of Technology, SE-971 87 Lule{\aa}, Sweden}
\email{\texttt{lech.maligranda@ltu.se}}
\maketitle


\section{Introduction and basic definitions}

For a Banach ideal space $X$ on $I = [0, 1]$ or $I = [0, \infty)$ let us consider, as in \cite{LM14}, the abstract Ces\`aro space 
$CX$ on $I$ defined as $CX = \{f\in L^0(I): C|f| \in X\}$ with the norm given by 
$$
\|f\|_{CX}=\|C|f| \|_{X},
$$
where $C$ is the Ces\`aro operator 
$$
Cf(x) = \frac{1}{x} \int_0^x f(t) \,dt, ~ x \in I.
$$
One may look at this space, on the one hand, as on generalizations of the well-known Ces\`aro spaces $Ces_p[0,1]$ and $Ces_p[0,\infty)$ 
which  were investigated for example in \cite{AM09} and on the other hand, just by definition, we get $C: CX\rightarrow X$ is bounded and $CX$ 
is the largest ideal space satisfying this relation, i.e. $CX$ is the optimal domain of $C$ for $X$. Consequently, the abstract Ces\`aro spaces 
may be considered also from the optimal domain point of view, as it was done in \cite{DS07}, \cite{NP10}, \cite{NP11}, \cite{MNS13}.
In this paper we discuss the Ces\`aro function spaces on $[0,\infty)$ and on $[0,1]$ from the point of view of optimal domain and 
optimal range of the Ces\`aro operator $C$. Such concept was already considered for $X = L^{p(\cdot)}$ on $[0,1]$ in \cite{NP10},  \cite{NP11} 
and for $X = L^{p(\cdot)}$ on $\mathbb{R}^n$ in \cite{MNS13}, althought the most interesting situation of $CX$ on  $[0,1]$ 
was omitted there. We develope and complete the discussion under some minimal assumptions. In this more interesting case of interval $[0,1]$ a very important role is played by the improvement of Hardy inequality presented in Theorem \ref{betterHardy}. 

We present some basic definitions to understand further description of results. 
A Banach space $X\subset L^{0} = L^{0}(I)$ is called a Banach ideal space on $I$ if $g\in X, f\in L^{0}(I), |f|\leq |g|$ a.e. on $I$ implies 
$f\in X$ and $|| f ||\leq ||g||$. We will also assume that ${\rm supp} X = I$, i.e. there exists $f \in X$ with $f(x) > 0$ for
each $x\in I$.  

For a given Banach ideal space $X$ on $I$ and a function $w\in L^0(I)$ such that $w(x) > 0$ a.e. on $I$, the weighted ideal 
Banach space $X(w)$ is defined as $X(w) = \{f\in L^0(I): fw \in X\}$ with the norm 
$$
\|f\|_{X(w)} =\| f w \|_{X}.
$$
In the whole paper only two concrete weights on $I = [0, 1]$ will appear, namely $v$ and $1/v$ where  $v(x)=1-x$. 
We will need also a non-increasing majorant $\widetilde{f}$ of a given function $f$, which is just 
$$
\widetilde{f}(x) = \esssup_{t \in I, \, t \geq x} |f(t)|,\ x\in I.
$$

Moreover, for a given Banach ideal space $X$ on $I$, we define a new Banach ideal space $\widetilde{X} = \widetilde{X}(I)$ as 
$\widetilde{X} = \{f\in L^0(I): \widetilde{f}\in X\}$ with the norm given by 
$$
\|f\|_{\widetilde{X}}=\|\widetilde{f}\|_{X}.
$$
By a {\it symmetric function space} on $I$ with the Lebesgue measure $m$ (symmetric space in short), we mean a Banach 
ideal space $X=(X,\| \cdot \|_{X})$ with the additional property that for any two equimeasurable functions 
$f \sim g, f, g \in L^{0}(I)$ (that is, they have the same distribution functions $d_{f}\equiv d_{g}$, where 
$d_{f}(\lambda) = m(\{x \in I: |f(x)|>\lambda \}),\lambda \geq 0,$ and $f\in E$ we have $g\in E$ and $\| f\|_{E} = \| g\|_{E}$. 
In particular, 
$\| f\|_{X}=\| f^{\ast }\|_{X}$, where $f^{\ast }(t)=\mathrm{\inf } \{\lambda >0\colon \ d_{f}(\lambda ) < t\},\ t\geq 0$. 

The {\it dilation operators} $\sigma_a$ ($a > 0$) defined on $L^0(I)$ by 
$$
\sigma_af(x) = f(x/a) \chi_{I}(x/a) = f(x/a) \chi_{[0, \, \min(1, \,a)]}(x), ~~ x \in I,
$$
are bounded in any symmetric space $X$ on $I$ and $\| \sigma_a\|_{X \rightarrow X} \leq \max (1, a)$ (see \cite[p. 148]{BS88} and 
\cite[pp. 96-98]{KPS82}). They are also bounded in some Banach ideal spaces which are not necessary symmetric spaces.
Furthermore, recall that the Ces\`aro operator $C$, the Copson operator $C^*$ and the Hardy-Littlewood maximal operator $M$ are defined, respectively, by 
$$
Cf(x) = \frac{1}{x}\int_0^x f(t) dt,\  x\in I, ~~ C^*f(x) = \int_{I\cap [x, \infty)} \frac{f(t)}{t} dt,\  x\in I,
$$ 
$$
Mf(x)=\sup_{a,b\in I,0\leq a\leq x\leq b}\frac{1}{b-a}\int_a^b |f(t)| dt,\  x\in I. 
$$

We refer the reader to \cite{LM14}, where basic facts about the spaces $CX$ and $\widetilde{X}$ were presented with more details. For more references on Banach ideal spaces and symmetric spaces we refer to \cite{KPS82}, \cite{LT79}, \cite{BS88}, \cite{KA77} and  \cite{Ma89}. 

\section{Optimal range}

Let $X$ and $Y$ be two Banach ideal spaces on $I$ and let $T: X\rightarrow Y$ be a bounded linear or sublinear operator. A Banach ideal 
space $Z$ on $I$ is called the {\it optimal domain} of $T$ for $Y$ within the class of Banach ideal spaces on $I$, if 
$T: Z\rightarrow Y$ is bounded and for each Banach ideal space $W$ on $I$, $T: W\rightarrow Y$ is bounded implies that 
$W\subset Z$. The last implication may be formulated equivalently as: if $Z$ and $W$ are Banach ideal spaces on $I$ and if 
$Z\subsetneq W$, then  $T: W\not \rightarrow Y$. Of course in such a case $X\subset W$. 

Similarly, we shall say that a Banach ideal space $Z$ on $I$ is the {\it optimal range} of $T$ for $X$ within the class of Banach ideal 
spaces on $I$, if $T: X\rightarrow Z$ is bounded and for each Banach ideal space $W$ on $I$, $T: X\rightarrow W$ is bounded implies 
that $Z\subset W$. Once again, the last condition may be replaced by:  $W \subsetneq Z$ implies $T: X\not\rightarrow W$. Such 
the optimal range satisfies of course $Z\subset Y$. 

The following theorem describes the optimal domain and optimal range problem for Ces\`aro operator within the class of Banach ideal spaces on $I$. 

\begin{theorem}\label{thm1}
Let $X$ be a Banach ideal space on $I$ such that the maximal operator $M$ is bounded on $X$.
\begin{itemize}
\item[(i)] If $I=[0,\infty)$, then $C: CX\rightarrow \widetilde{X}$ is bounded.
Moreover, the space $CX$ is the optimal domain of $C$ for $X$ and for $\widetilde{X}$ (also for $CX$ if the dilation operator $\sigma_a$ is bounded on $X$ for some $0 < a < 1$). The space $\widetilde{X}$ is the optimal range of $C$ for $CX$, $X$ and $\widetilde{X}$. In particular, $CX=C\widetilde{X}$. 
\item[(ii)] If $I = [0,1]$, then $C: CX\rightarrow \widetilde{X(1/v)}(v)$ 
is bounded. The space $CX$ is the optimal domain of $C$ for $X$ and also for $\widetilde{X(1/v)}(v)$.
Moreover, if the maximal operator $M$ is bounded on $X'$, then the space 
$\widetilde{X(1/v)}(v)$ is the optimal range of $C$ for $CX$ and $X(v)$ {\rm (cf. Diagram 2)}. In particular, 
$CX=C[\widetilde{X(1/v)}(v)]$.
 
\item[(iii)] If $I=[0,1]$ and the dilation operator $\sigma_{1/2}$ is bounded on $X$, then $C: C\widetilde{X}\rightarrow \widetilde{X}$ is bounded. Moreover, the space $C\widetilde{X}$ is the optimal domain of $C$ for $\widetilde{X}$ and the space $\widetilde{X}$ is 
the optimal range of $C$ for $C\widetilde{X}$, $X$ and $\widetilde{X}$. One also has $C\widetilde{X}=CX\cap L^1$. 
\end{itemize} 
\end{theorem}

Before we prove the theorem, let us comment the situation. Suppose that the corresponding assumptions in Theorem \ref{thm1} are satisfied. Of course, boundedness of $M$ on $X$ implies also boundedness of $C$ on $X$, therefore support of $CX$ is for sure the same as support of $X$ (cf. \cite{LM14}).   
Let $I=[0,\infty)$. Then the statement of (i) may be therefore pictured, putting the boundedness of $C$ and respective embeddings, on the following diagram. 
\vspace{3mm}

\[
\begin{tikzcd}
CX\arrow{r}{C} 
&\widetilde{X}\arrow[hook]{r}
&X\arrow[hook]{r}
&CX\\
X\arrow[hook]{u}  \arrow{ru}{C}\\
\widetilde{X} \arrow[hook]{u} \arrow{ruu}{C}
\end{tikzcd}
\]
\centerline {Diagram 1}
\vspace{3mm}

Moreover, point (i) says that, in fact, $CX$ is the optimal domain of $C$ for $\widetilde{X}$, since  $CX=C\widetilde{X}$. 
Even more can be said when the dilation operator $\sigma_a$ is bounded on $X$ for a certain $0<a<1$. Then  $CX$ is 
the optimal domain of $C$ even for $CX$ since, by Lemma 6 in \cite{LM14}, it follows that $CCX = CX$. On the other 
hand, we will see that $\widetilde{X}$ is the  optimal range of $C$ for $\widetilde{X}$, which by the above diagram 
means that also for $X$ and for $CX$. 

Much more interesting and delicate is the case of interval $[0,1]$.
Suppose that $C: X\rightarrow X$ is bounded and all assumptions of (ii) and (iii) are satisfied. Then $C: CX\rightarrow X$ is bounded, 
where $CX$ is by definition the optimal domain of $C$ for $X$. The case (ii) says that the optimal range of $C$ for $CX$ is then 
$\widetilde{X(1/v)}(v)$. It is however interesting that one may look at the situation also in another way. 
Let's start once again with $C:X\rightarrow X$ and find first the optimal range. It appears to be just $\widetilde{X}$ 
(cf. \cite[Theorem 8.2]{NP10}, \cite[Theorem 3.16]{NP11} and \cite[Theorem 4.1]{MNS13}) which is much smaller than $\widetilde{X(1/v)}(v)$. If we now 
find optimal domain of $C$ for $\widetilde{X}$ it is then just $CX\cap L^1 = C(\widetilde{X})$. 
The diagram describing this dichotomy is now more complicated. 
\vspace{3mm}

\[
\begin{tikzcd}
X(v)\arrow[hook]{r}
&CX\arrow{r}{C} 
& \widetilde{X(1/v)}(v) \arrow[hook]{r}
&X\arrow[hook]{r}
& C\widetilde{X}\arrow[hook]{r}
&CX\\
X \arrow[hook]{r} \arrow[hook]{ru}\arrow[hook]{u}
&C\widetilde{X} \arrow[hook]{u} \arrow{r}{C} \arrow{ru}{C} 
&\widetilde{X}\arrow[hook]{u}\\
\widetilde{X(1/v)}(v) \arrow[hook]{u}\arrow{rru}{C} \\
\widetilde{X}\arrow[hook]{u}\arrow{rruu}{C} 
\end{tikzcd}
\]
\centerline {Diagram 2}
\vspace{3mm}

In general, there is no inclusion relation between $X(v)$ and $C\widetilde{X}$. For example, if $X$ is a symmetric space on 
$I = [0, 1]$, we have for $f(x): = \frac{1}{1-x}$ that $f\in X(v)$ while $f\not \in C\widetilde{X}$ because $Cf(x)\rightarrow \infty$ as 
$x \rightarrow 1^-$ and so $\widetilde {Cf}$ is not defined (or just $\infty$ everywhere). Therefore, $X(v)\not\subset C\widetilde{X}$. 
This means also that $C$ does not act from $X(v)$ into $\widetilde{X}$. On the other hand, let $X=L^2$ and put 
$f(x)=|\frac{1}{2}-x |^{-1/2}$. Then $f\not\in L^2$, but $Cf\in L^{\infty}$ and so $\widetilde{Cf}\in L^{\infty}\subset L^2$. This gives 
$ C\widetilde{X}\not\subset X(v)$. For general symmetric space $X$ on $I$ such that $C: X\rightarrow X$ is bounded, one could 
take $f\in L^1$ in such a way that $f-f\chi_{[1/2-\epsilon,1/2+\epsilon]}\in L^{\infty}$ for each $0<\epsilon <1/2$ but $f\not\in X$, to achive the same effect. 

\proof[Proof of Theorem 1]
(ii). Let $0\leq f\in CX$. Suppose first that $0\leq y\leq t\leq 2y\leq 1$. Then 
\begin{equation}
Cf(t)= \frac{1}{t}\int_0^t f(s) ds\geq \frac{1}{2y}\int_0^y f(s) ds = \frac{1}{2}Cf(y). \label{R1}
\end{equation}
If now $0 \leq x\leq y$ and $y\leq \frac{1}{2}$, then applying (\ref{R1}) one gets 
\begin{equation*}
MCf(x)\geq \frac{1}{2y-x}\int_x^{2y} Cf(t)dt\geq \frac{1}{2y}\int_y^{2y} Cf(t)dt \label{R2a}
\end{equation*}
\begin{equation*}
\geq \frac{1}{2y}\int_y^{2y} \frac{Cf(y)}{2}dt=\frac{1}{4}Cf(y)\geq \frac{1-y}{4(1-x)}Cf(y). \label{R2}
\end{equation*}
Suppose now that $\frac{1}{2}\leq y\leq t\leq 1$. Then, similarly as in (\ref{R1}),
\begin{equation}
Cf(t)= \frac{1}{t}\int_0^t f(s) ds\geq \int_0^y f(s) ds \geq \frac{1}{2}Cf(y). \label{R3}
\end{equation}
In consequence, when $0 \leq x\leq y$ and $ \frac{1}{2}\leq y \leq 1$, applying (\ref{R3}) we obtain
\begin{equation*}
MCf(x)\geq \frac{1}{1-x}\int_x^1 Cf(t)dt\geq \frac{1}{1-x}\int_y^1 Cf(t)dt
\end{equation*}
\begin{equation*}
\geq \frac{1}{1-x}\int_y^1 \frac{Cf(y)}{2}dt= \frac{1-y}{2(1-x)}Cf(y).
\end{equation*}
Consequently, 
\begin{equation}
MCf(x)\geq \frac{1}{4(1-x)} \esssup_{0 \leq x\leq y \leq 1} (1-y)Cf(y)=\frac{1}{4(1-x)} \widetilde{[vCf]}(x).
\end{equation}
Since $M$ is bounded on $X$, by our assumption, it follows that
\begin{equation*}
\|Cf\|_{\widetilde{X(1/v)}(v)}= \|\widetilde{[vCf]}/v\|_X\leq 4\|M\|_{X\rightarrow X}\|Cf\|_X=4\|M\|_{X\rightarrow X}\|f\|_{CX}.
\end{equation*}
This means that $C:CX\rightarrow \widetilde{X(1/v)}(v)$ is bounded and the first statement of (ii) is proved. 
It remains to show that the space $\widetilde{X(1/v)}(v)$ is optimal range of $C$ for $CX$ (in fact, even for $X(v)$). Suppose that there is a Banach ideal space $Z$ on $I$ such that 
$$
Z\subsetneq Y ~{\rm but}~ C: CX\rightarrow Z ~{\rm is ~bounded}.
$$
Let $ 0\leq f\in Y \backslash Z$. Define 
$$
g(x)=\frac{1}{(1-x)} \widetilde{[vf]}(x), x \in I.
$$
Then $f \leq g$ and $g\in \widetilde{X(1/v)}(v)\subset X$ because $\frac{1}{1-x} \widetilde{[vg]}(x) = \frac{1}{1-x} \widetilde{[vf]}(x)$. We have 
\begin{eqnarray*}
C(g/v)(x)
&=&
\frac{1}{x}\int_0^{x} \frac{\widetilde{[vg]}(t)}{(1-t)^2}dt \geq \frac{\widetilde{[vf]}(x)}{x}\int_0^{x} \frac{1}{(1-t)^2}dt \\
&=&
\frac{\widetilde{[vf]}(x)}{x}\frac{x}{(1-x)}\geq f(x),
\end{eqnarray*}
which means that $C(g/v)\not \in Z$. However, $g\in X$ and so $g/v\in X(v)$. Also, by Theorem \ref{betterHardy}, 
$X(v)\subset CX$ and therefore $g/v \in CX$ which means that $C: CX\not \rightarrow Z$. Note that we have already 
shown $C: X(v)\not \rightarrow Z$, which by inclusion  $X(v)\subset CX$ means that $\widetilde{X(1/v)}(v)$ is the optimal range also for $X(v)$. 

(iii). The argument is analogous to the one from statement (5.1) in \cite{NP10}. However, we need to modify it because in \cite{NP10} the maximal operator is defined on a larger interval than $[0,1]$. 
Let $0\leq f \in CX\cap L^1[0,1]$. We shall understand that $f(x) = 0$ for $x > 1$. Of course, inequality from (\ref{R1}) remains true in this case, since  $f\in L^1[0,1]$. Suppose  that $0<x \leq y \leq 1$ and consider two cases. If $y/2\leq x$, then 
$$
M\sigma_ {1/2} Cf(x)\geq \frac{2}{y}\int_{y/2}^{y} \sigma_{1/2} Cf(u)du.
$$
If $x\leq y/2$, then 
$$
M\sigma_{1/2} Cf(x)\geq \frac{1}{y-x}\int_{x}^{y} \sigma_{1/2} Cf(u)du\geq \frac{1}{y}\int_{y/2}^{y} \sigma_{1/2} Cf(u)du.
$$
Alltogether we get 
$$
M\sigma_{1/2} Cf(x) \geq \frac{1}{y}\int_{y/2}^{y} \sigma_{1/2} Cf(u)du = \frac{1}{2y}\int_{y}^{2y} Cf(t)dt\geq \frac{1}{4}Cf(y).
$$
Therefore, similarly as before, 
\begin{equation*}
M\sigma_{1/2} Cf(x)\geq \frac{1}{4}  \esssup_{x\leq y}Cf(y)= \frac{1}{4} \, \widetilde{Cf}(x),
\end{equation*}
which gives 
\begin{eqnarray*}
\|f\|_{C\widetilde{X}} 
&=& 
\|\widetilde{Cf}\|_X\leq 4 \, \|M\sigma_{1/2} Cf\|_X\leq 4 \, \|M\|_{X\rightarrow X}\|\sigma_{1/2} \|_{X\rightarrow X}\|Cf\|_X \\
&=& 
4 \, \|M\|_{X\rightarrow X}\|\sigma_{1/2} \|_{X\rightarrow X}\|f\|_{CX} \leq 4 \, \|M\|_{X\rightarrow X}\|\sigma_{1/2} \|_{X\rightarrow X}\|f\|_{CX\cap L^1}.
\end{eqnarray*}
On the other hand, if $0\leq f \in C\widetilde{X}$, then
$$
\|f\|_{L^1} =\int_0^1f(t) dt\frac{\|\chi_{[0,1]}\|_X}{\|\chi_{[0,1]}\|_X}= \frac{\|(\int_0^1f(t) dt)\chi_{[0,1]}\|_X}{\|\chi_{[0,1]}\|_X}
\leq \frac{\|\widetilde{Cf}\|_X}{\|\chi_{[0,1]}\|_X}.
$$
Thus also
$$
\|f\|_{CX\cap L^1}\leq \max\{1, \frac{1}{\|\chi_{[0,1]}\|_X}\}\|\widetilde{Cf}\|_X,
$$
which means that $C\widetilde{X}=CX\cap L^1$. 
For the sake of completeness we present the argument that $\widetilde{X}$ is the optimal range of $C$ for $C\widetilde{X}$, although it 
works just like in \cite[Theorem 8.2]{NP10}.  Let $Z$ be a Banach ideal space on $I$ and suppose that $0\leq f\in \widetilde{X} \backslash Z$. 
Then also $\widetilde{f}\in \widetilde{X} \backslash Z$ and $C\widetilde{f}\geq \widetilde{f}$. However $\widetilde{f}\not \in Z$, which 
means that  $C\widetilde{f}\not \in Z$ and  $C: C\widetilde{X}\not \rightarrow Z$. 

(i) This case is easier and may be deduced directly from \cite{MNS13}. Since for 
$0<y$ also $2y\in I$ it is enough to follow (\ref{R1}) and after that to get  for $y \geq x \geq 0$
\begin{equation*}
MCf(x)\geq \frac{1}{2y-x}\int_x^{2y} Cf(t)dt\geq \frac{1}{4}Cf(y).
\end{equation*}
Then
\begin{equation*}
\|Cf\|_{\widetilde{X}} =\|\widetilde{Cf}\|_X\leq 4\|MCf\|_X\leq 4\|M\|_{X\rightarrow X}\|Cf\|_X=4\|M\|_{X\rightarrow X}\|f\|_{CX},
\end{equation*}
which means that $C:CX\rightarrow \widetilde{X}$ is bounded and $CX = C\widetilde{X}$. The optimal range of $C$ for 
$\widetilde{X},X,CX$ is once again $\widetilde{X}$ and the proof is the same as in (iii) (see also \cite[Theorem 8.2]{NP10}, 
\cite[Theorem 3.16]{NP11} and \cite[Theorem 4.1]{MNS13}).
\endproof

\section{Hardy inequality} 
 We present an improvement of the Hardy inequality which appear for spaces on $I=[0,1]$.
\begin{theorem}\label{betterHardy}
If $C$ is bounded on a Banach ideal space $X$ on $I=[0,1]$ and maximal operator $M$ is bounded on $X'$, then 
\[
C:X(v)\rightarrow X
\]
is also bounded.
\end{theorem}

\proof Let $0\leq f\in X$. We have for $0< x\leq \frac{1}{2}$ 
\[
C(f/v)(x)=\frac{1}{x}\int_0^x\frac{f(s)}{1-s}ds\leq \frac{2}{x}\int_0^xf(s)ds
\]
and for $ \frac{1}{2} < x\leq 1$
\[
C(f/v)(x)=\frac{1}{x}\int_0^x\frac{f(s)}{1-s}ds\leq 2\int_0^x\frac{f(s)}{1-s}ds.
\]
If we define an operator $T$ as $Tf(x)=\int_0^x\frac{f(s)}{1-s}ds$, then
\[
C(f/v)\leq 2(Cf + Tf).
\]
Therefore, we need to show that $T$ is bounded on $X$. 
Consider an involution operator $\tau:f(x)\mapsto f(1-x)$. Then 
\begin{equation}
Tf(x)=\int_0^x\frac{f(s)}{1-s}ds=\int_{1-x}^1\frac{f(1-s)}{s}ds=\tau C^* \tau f(x). 
\end{equation}
Observe that the space 
\[
X^-=\{f:\tau f\in X\} 
\]
with its natural norm $\|f\|_{X^-}=\|\tau f\|_{X}$ is also a Banach ideal space on $I$ and $(X^-)^-$. Just by definition $\sigma :X\rightarrow X^-$,  $\tau :X^-\rightarrow X$ are bounded and $\tau \tau = id$. Thus $T$ is bounded on $X$ if and only if $C^*$ is bounded on $X^-$. We will prove the last equivalence. 
Notice that simply 
\begin{align}
Mf(1-x)&=\sup_{a\not =b,0\leq a\leq 1-x\leq b\leq 1}\frac{1}{b-a}\int_a^bf(s)ds\\&=\sup_{a\not =b,0\leq 1-b\leq x\leq 1-a\leq 1}\frac{1}{b-a}\int_{1-b}^{1-a}f(1-s)ds=(M\tau f)(x)
\end{align}
and so $M\tau f=\tau Mf$ which means that for any Banach ideal space $Y$, $M$ is bounded on $Y$ if and only if M is bounded on $Y^-$, which by our assumption gives that $M$ is bounded on $(X')^-$. Thus also $C$ is bounded on $(X')^-$ and by duality $C^*$ is bounded on $[(X')^-]'$. However, it is evident that for any Banach ideal space $Y$ there holds  $(Y')^-=(Y^-)'$. Then $[(X')^-]'=(X'')^-=X^-$ and so $C^*$ is bounded on $X^-$. 
\endproof
\begin{remark}
If $X$ is a symmetric space, then evidently $X=X^-$ and we get Lemma 10 from \cite{LM14}, which proof was a generalization of the Astashkin - Maligranda result from \cite{AM09}. Moreover, our Theorem 2 includes Theorem 9 in \cite{LM14} for the weighted $L^p(x^{\alpha})$ spaces when $1\leq p< \infty$ and $-1/p<\alpha <1-1/p$.    
\end{remark}


\begin{thebibliography}{99}

\bibitem[AM09]{AM09} S.V. Astashkin and L. Maligranda, {\it Structure of Ces\'aro function spaces}, Indag. Math. (N.S.) 20 (2009), 
no. 3, 329--379.

\bibitem[BS88]{BS88} C. Bennett and R. Sharpley, {\it Interpolation of Operators,} Academic Press, Boston 1988.

\bibitem[DS07]{DS07} O. Delgado and J. Soria, {\it Optimal domain for the Hardy operator}, J. Funct. Anal. 244 (2007), no. 1, 119--133.

\bibitem[KA77]{KA77} L. V. Kantorovich and G. P. Akilov, {\it Functional Analysis}, Nauka, Moscow 1977 (Russian); English 
transl. Pergamon Press, Oxford-Elmsford, New York 1982.

\bibitem[KMS07]{KMS07} R. Kerman, M. Milman and G. Sinnamon, {\it On the Brudnyi-Krugljak duality theory of spaces formed 
by the K-method of interpolation}, Rev. Mat. Complut. 20 (2007), no. 2, 367--389.

\bibitem[KPS82]{KPS82} S. G. Krein, Yu. I. Petunin, and E. M. Semenov, {\it Interpolation of Linear Operators}, Amer. Math. Soc., 
Providence 1982.

\bibitem[LM14]{LM14} K. Le\'snik and L. Maligranda, {\it On abstract Ces\`aro spaces. I. Duality}, preprint of 20 pages, 20 March 2014, {\tt arXiv:1401.6415v2} at: http://arxiv.org/pdf/1401.6415.pdf.

\bibitem[LT79]{LT79} J. Lindenstrauss and L. Tzafriri, {\it Classical Banach
Spaces, II. Function Spaces}, Springer-Verlag, Berlin-New York 1979.


\bibitem[Ma89]{Ma89} L. Maligranda, {\it Orlicz Spaces and Interpolation}, Seminars 
in Mathematics 5, University of Campinas, Campinas SP, Brazil 1989. 

\bibitem[MNS13]{MNS13} Y. Mizuta, A. Nekvinda and T. Shimomura, {\it Hardy averaging operator on generalized Banach function spaces and duality}, Z. Anal. Anwend. 32 (2013), no 2, 233--255.

\bibitem[NP10]{NP10} A. Nekvinda and L. Pick, {\it Optimal estimates for the Hardy averaging operator}, Math. Nachr. 283 (2011), no. 2, 262--271.

\bibitem[NP11]{NP11} A. Nekvinda and L. Pick, {\it Duals of optimal spaces for the Hardy averaging operator}, Z. Anal. Anwend. 30 (2011), no. 4, 435--456.


\end{thebibliography}
\end{document}